\documentclass[10pt,reqno]{amsart}
\usepackage[letterpaper,textwidth=6.12in,textheight=8.94in,centering]{geometry}
\usepackage[T1]{fontenc}
\usepackage{lmodern}
\usepackage{amsmath,amssymb}
\usepackage{graphicx,booktabs,float,enumitem,placeins}
\usepackage[numbers,sort&compress]{natbib}
\usepackage[hidelinks]{hyperref}
\theoremstyle{plain}
\newtheorem{theorem}{Theorem}[section]
\newtheorem{lemma}[theorem]{Lemma}

\newtheorem{proposition}[theorem]{Proposition}
\theoremstyle{definition}
\newtheorem{definition}[theorem]{Definition}

\newtheorem{example}[theorem]{Example}

\theoremstyle{remark}
\newtheorem{remark}[theorem]{Remark}

\newenvironment{algorithmindent}
  {\begin{list}{}{%
     \setlength{\leftmargin}{1.3em}%
     \setlength{\topsep}{1pt}%
     \setlength{\parsep}{0pt}%
     \setlength{\listparindent}{0pt}%
     \setlength{\partopsep}{0pt}%
   }\item\relax}
  {\end{list}}

\setlist[enumerate,1]{label=(\roman*),leftmargin=*}
\floatstyle{ruled}
\newfloat{algorithm}{tbp}{loa}
\floatname{algorithm}{Algorithm}
\title[Decomposition algorithms for rational curves]{Decomposition Algorithms for Rational Curves on Real Classical Groups}
\author{Yifan Li}
\address{State Key Laboratory of Mathematical Sciences, AMSS, University of Chinese Academy of Sciences, Beijing, 100190, China}
\email{liyifan@amss.ac.cn}

\author{Zijia Li}
\address{State Key Laboratory of Mathematical Sciences, AMSS, University of Chinese Academy of Sciences, Beijing, 100190, China}
\email{lizijia@amss.ac.cn}

\author{Ke Ye}
\address{State Key Laboratory of Mathematical Sciences, AMSS, University of Chinese Academy of Sciences, Beijing, 100190, China}
\email{keyk@amss.ac.cn}

\date{}
\thanks{\copyright~2026. Copyright held by the owner/author(s).}
\keywords{Involution, Multiplicity, Orthogonal matrices, Fundamental theorem, Linear/Quadratic equations, Degree and Poles}

\begin{document}

\begin{abstract}
This paper is devoted to an algorithmic investigation of the decomposition problem for rational curves on real classical groups. The main contributions are: 
    (1) We present two Maple-implemented\footnotemark{} 
    algorithms for decomposing rational curves on real classical groups. The first algorithm decomposes a curve with multiple poles into a unique product of curves with single poles. For rational curves on (indefinite) orthogonal groups, the second algorithm further decomposes single-pole curves into products of low-degree curves. The former only involves linear equations, whereas the latter requires solving both linear and quadratic equations. 
    (2) We compare the decomposition of the same curve with respect to different ambient groups. Through this comparison, we provide an interpretation of the "multiplication trick" for motion polynomials within our framework. Moreover, we demonstrate that our second algorithm aligns with this trick. 
    (3) We prove that the solvability of linear equations in the second algorithm is preserved under the extension of the ambient group. 
For illustration, we provide comprehensive computational examples.
\end{abstract}
\maketitle
\setcounter{footnote}{1}
\footnotetext{The Maple codes are available at
  \url{https://github.com/Joshua-1417/decomposing-rational-curves/}}


\section{Introduction}\label{sec:introduction}
Rational curves are fundamental objects of study across both pure and applied mathematics. Within pure mathematics, they are indispensable in modern algebraic geometry, serving as essential tools in areas such as the minimal model program~\cite{Kollar92,Mori1988flip}, the study of rational and unirational varieties~\cite{Graber03families}, and the classification of Fano varieties~\cite{Kollar92fano}. They also play a central role in a variety of enumerative problems in real algebraic geometry~\cite{KR15,KR17,Mikhalkin17,Sottile00,Welschinger05}. Beyond theory, rational curves have wide-ranging applications in engineering practice~\cite{Roth79tk,Darboux1890}. They provide a classical framework in kinematics for parameterizing and analyzing motions, and they are imperative in computer-aided geometric design for the efficient modeling of complex shapes and objects~\cite{Farin97, Juttler93zwanglaufige, Winkler08}.

The decomposition problem for rational curves asks whether every rational curve can be decomposed into a product of low-degree curves (cf. Definition~\ref{def:deg}). This problem dates back to Kempe in 1875 \cite{MR1575631}. The past three decades have witnessed a revival of this problem, driven by advances in techniques and methods. For instance, Kempe's result was reformulated and greatly generalized by a topological method \cite{KM96,KM98,KM02}; the motion polynomial was introduced to study the decomposition problem using techniques from geometric algebra \cite{Hegedus13, LSS19, GKLRSV17,LJS18,Thimm2026}; a matrix model for this problem was proposed, extending Kempe's result to homogeneous spaces \cite{Liye2024rational}.

\subsection*{Main results}
In contrast to the purely theoretical study in \cite{Liye2024rational} of the decomposition problem for rational curves, this paper focuses on algorithmic aspects. We present two decomposition algorithms: the first (cf. Algorithm~\ref{alg:decomposition1}) decomposes a curve with multiple poles into a product of curves with single poles; the second (cf. Algorithm~\ref{alg:decomposition2}) further decomposes curves on (indefinite) orthogonal groups with single poles into a product of low-degree curves. The correctness of both Algorithms~\ref{alg:decomposition1} and \ref{alg:decomposition2} is theoretically ensured (cf. Theorems~\ref{alg1:correctness} and~\ref{alg2:correctness}). They are efficiently implemented in Maple. For illustration, we provide several examples throughout this paper.

Furthermore, we prove that the decomposition obtain by Algorithm~\ref{alg:decomposition1} is unique up to a permutation of factors (cf. Proposition~\ref{uniq}). By comparing the decomposition of a rational curve in different ambient groups, we obtain an interpretation (cf. Subsection~\ref{subsec:ext}) of the mysterious "multiplication trick" in the method of motion polynomials \cite{GKLRSV17} within our framework. We also show that the solvability of linear equations in Algorithm~\ref{alg:decomposition2} is invariant under extension of the ambient group (cf. Theorem~\ref{nece}).

\subsection*{Organization of the paper}
Section~\ref{sec:pre} introduces basic notations and results. In Section~\ref{sec:decomp1}, we present Algorithm~\ref{alg:decomposition1}, validate its correctness, and illustrate it with examples. In Section~\ref{sec:decomp2}, we present Algorithm~\ref{alg:decomposition2} and discuss some of its properties. In Section~\ref{sec:examples}, we provide additional examples. Section~\ref{sec:summary} provides a brief summary and further discussion.

\section{Preliminary} \label{sec:pre}
We begin by recalling the definition of rational curves on a real algebraic variety. These curves will be the central objects of our study.
\begin{definition} [Rational Curve]
    Let $X \subseteq \mathbb{R}^m $ be a real variety, a rational curve on $X$ is a rational map $  \gamma : \mathbb{P}^1 \rightarrow X$.
\end{definition}
Concretely, a rational curve $\gamma$ on $X$ can be expressed in affine coordinates as
\begin{equation*}
   \begin{aligned}
    &\gamma : \mathbb{R}^1 \rightarrow X, \, \gamma(t) = \left( \frac{p_1(t)}{q(t)}, \cdots , \frac{p_m(t)}{q(t)} \right),\\[5pt]
   & \gamma(+\infty) = \gamma(- \infty)  \in X. 
    \end{aligned} 
\end{equation*}
where $p_1, \cdots p_m, q \in \mathbb{R}[t]$, $\text{gcd}(p_1, \cdots , p_m, q) = 1$, and $q(t)$ has no real roots, ensuring the map is well-defined over $\mathbb{R}$.
We now introduce two fundamental invariants of rational curve.
\begin{definition} [Degree and Pole]\label{def:deg}
The \textbf{degree} of a rational curve $\gamma$ on $X$ is defined as $\text{deg}(q)$. If $\zeta \in \mathbb{C}$ such that $q(\zeta) = 0$, then called $\zeta$ is a \textbf{pole} of $\gamma$.
\end{definition}
In this work, we focus on rational curves defined on a specific and important class of real algebraic groups, described as follows
\begin{equation*}
    G_B(\mathbb{F}) = \{ X \in \mathbb{F}^{n \times n} : X B X^{\sigma} = B \}.
\end{equation*}
where $\mathbb{F} = \mathbb{R}, \mathbb{C}$ or $\mathbb{H}$, $\sigma$ denotes transpose or conjugate transpose, and $B  \in \operatorname{GL}_n(\mathbb{F})$ is a fixed invertible matrix. This framework encompasses several classical groups, as illustrated by the following examples.
\begin{enumerate}[label= (\arabic*)]
    \item $\mathbb{F} = \mathbb{R}$, $\sigma = $ transpose, $B = I_n$, $G_B(\mathbb{F}) = O_n(\mathbb{R})$.
    \item $\mathbb{F} = \mathbb{C}$, $\sigma = $ conjugate transpose, $B = I_n$, $G_B(\mathbb{F}) = U_n$.
    \item $\mathbb{F} = \mathbb{R}$, $\sigma = $ transpose, $B =  \text{diag}(I_p, -I_{n-p})$, $G_B(\mathbb{F}) = O_{p,n-p}(\mathbb{R})$.
\end{enumerate}
We will study rational curves $\gamma : \mathbb{P}^1 \rightarrow G_B(\mathbb{R})$ that are normalized to satisfy $\gamma(\infty) = I_n$. Such a curve can be expressed in matrix form as
\begin{equation*}
\frac{\begin{pmatrix}
p_{11}(t) & p_{12}(t) & \cdots & p_{1n}(t) \\[4pt]
p_{21}(t) & p_{22}(t) &  \cdots & p_{2n}(t) \\[4pt]
\vdots & \vdots & \ddots & \vdots \\[4pt]
p_{n1}(t) & p_{n2}(t) & \cdots & p_{nn}(t)
\end{pmatrix}}{q(t)} 
\end{equation*} 
with $\text{deg}(\gamma) = \text{deg}(q) = \text{deg}(p_{ii})>\text{deg}(p_{ij})$ for $i \neq j$, and $q$ has no real roots.
\begin{remark}
The condition that $q$ has no real roots implies that the degree of any rational curve on $G_B(\mathbb{R})$ is an even number.
\end{remark}
To make these definitions clear, consider the following example of a rational curve on $SO_3(\mathbb{R})$
\begin{equation*}
\gamma(t) = \frac{\begin{pmatrix}
t^{4} - 3t^{2} - 4 & 4t^{3} + 4t & 0 \\[3pt]
-4t^{3} + 4t & t^{4} - 5t^{2} + 4 & 2t^{3} + 8t \\[3pt]
8t^{2} & -2t^{3} + 8t & t^{4} + 3t^{2} - 4
\end{pmatrix}}{(t^2+1)(t^2+4)}.
\end{equation*}
One can verify that $\gamma(t)\gamma(t)^\mathsf{T} = I_3$. And $\text{deg}(\gamma) = 4$, its poles in $\mathbb{C}$ are $\{\pm i, \pm 2i \}$. 
We note a useful algebraic property of this set of rational curves. Let $\gamma$ be rational curve on $G_B(\mathbb{R})$ and  with $\gamma(\infty) = I_n$, we have $\gamma B \gamma^{\sigma} = B$, then $\gamma^{-1} = B \gamma^{\sigma} B^{-1}$, then all the rational curves on $G_B(\mathbb{R})$ with $\gamma(\infty) = I_n$ form a group under multiplication. 
The structure of this group of rational curves is greatly clarified by two key results. The first is a fundamental operation that allows us to reduce the degree of a curve by "removing" a pair of conjugate poles.
\begin{lemma} [Degree Reduction \cite{Liye2024rational}]
\label{deg_red}
For any rational curve $\gamma(t) = P(t)/q(t)$ on $G_B(\mathbb{R})$ of degree $d$ with a pole $\zeta \in \mathbb{C}/\mathbb{R}$ of multiplicity $s$, there exists a rational curve $\alpha(t)$ on $G_B(\mathbb{R})$ of degree $2s$ with only poles $c$ and $\overline{c}$ such that $\alpha(t)\gamma(t)$ is a rational curve on $G_B(\mathbb{R})$ of degree $d-2s$.
\end{lemma} 

Iterative application of the degree reduction principle leads to the following central decomposition theorem, which asserts that every rational curve can be factored into curves corresponding to its pairs of conjugate poles.

\begin{theorem} [Decomposition Theorem \cite{Liye2024rational}]
\label{decom_thm}
Suppose $\gamma$ is a rational curve on $G_B(\mathbb{R})$, which has poles at $\zeta_1, \overline{\zeta_1}, \cdots \zeta_s, \overline{\zeta_s}$ with multiplicity $d_1, d_1, \cdots ,d_s, d_s$, then 
\begin{equation*}
 \gamma = \prod_{j=1}^s \gamma_j.   
\end{equation*}
for rational curve $\gamma_j$ on $G_B(\mathbb{R})$ with poles at $\zeta_j, \overline{\zeta_j}$  of multiplicity $d_j, d_j$. 
  \end{theorem}

\begin{remark}
  If a rational curve $\lambda$ on $G_B(\mathbb{R})$ only has a pair of conjugate poles $\zeta$ and $\overline{\zeta}$, we say $\lambda$ is of single pole, and if rational curve $\gamma$ on $G_B(\mathbb{R})$ has more than one pair of conjugate poles, we say that $\gamma$ is of multiple poles. And in the following content, if $\gamma$ has a pair of conjugate poles $\zeta$ and $\overline{\zeta}$, then we will simply write $\gamma$ has pole $\zeta$ for simplicity.
\end{remark}

A natural question arises concerning the uniqueness of the decomposition provided by Theorem~\ref{decom_thm}. While the theorem guarantees existence, the uniqueness of the decomposition is not immediate. We establish this uniqueness under a natural condition—namely, the decomposition is unique once an order for the poles is fixed.

\begin{proposition}\label{uniq}
Let $\gamma(t) = P(t)/q(t)$ be a rational curve on $G_B(\mathbb{R})$, suppose $q(t) = q_1(t)\cdots q_s(t)$, where $q_i(t)$ is powers of irreducible quadratic polynomial over $\mathbb{R}$. If the order of $q_1(t), \cdots , q_s(t)$ is chosen, then the decomposition of $\gamma = \left(P_1(t)/q_1(t) \right) \cdots \left(P_s(t)/q_s(t)\right)$ is unique.
\end{proposition}

\begin{proof}
Suppose $\gamma(t) = \left(P_1(t)/q_1(t) \right) \cdots \left(P_s(t)/q_s(t)\right)$ and $\gamma(t) = \left(R_1(t)/q_1(t) \right) \cdots \left(R_s(t)/q_s(t)\right)$ are decompositions of $\gamma(t)$. We make an induction on $s$. If $s = 1$, the case is trivial, it is obviously true.
 Suppose the decomposition is unique for $s-1$, we will prove it for the case $s$. First we have
\begin{equation}\label{eq_for_uniq}
    \frac{R_1'(t)}{q_1(t)} \frac{P_1(t)}{q_1(t)} \cdots \frac{P_s(t)}{q_s(t)} = \frac{R_2(t)}{q_2(t)} \cdots \frac{R_s(t)}{q_s(t)}.
\end{equation}
where $R_1'(t)/q_1(t)$ is the inverse of rational curve $R_1(t)/q_1(t)$.
 Let $\gamma'$ be the rational curve on the right side of equation \eqref{eq_for_uniq}, then the poles of $\gamma'$ are the roots of  $q_2(t), \cdots , q_s(t)$. And for the rational curve $\alpha = \left(R_1'(t)P_1(t)\right)/q_1^2(t)$, the pole of it is the root of $q_1(t)$. Since the poles on both sides of the equation must remain the same, hence $q_1^2(t) \mid  R_1'(t)P_1(t)\cdots P_s(t)$. Let $P'(t) = P_2(t)\cdots P_s(t)$, we have
\begin{equation*}
   P'(t) B P'^{\sigma}(t) = q_2^2(t) \cdots q_s^2(t) B. 
\end{equation*}
Then $q_1^2(t) \nmid  P'(t)$, otherwise $q_1^2(t) \mid  P'(t) B P'^{\sigma}(t)$, contradiction. Thus $q_1^2(t) \mid  R_1'(t)P_1(t)$.
In fact, we have $q_1^2(t) I_n = R_1'(t)P_1(t)$, which implies that $P_1(t) = R_1(t)$. And consider
\begin{equation*}
   \gamma' =  \frac{R_1'(t)}{q_1(t)} \gamma = \frac{R_2(t)}{q_2(t)} \cdots \frac{R_s(t)}{q_s(t)}. 
\end{equation*}
which only have $s-1$ poles, then it can be decomposed uniquely by induction.  Consequently $P_i(t) = R_i(t)$, we show that the decomposition is unique for $s$ poles.
\end{proof}

\section{Decomposition Algorithm for Rational Curves with Multiple Poles}\label{sec:decomp1}
Let $\gamma = P(t)/q(t)$ be a rational curve on $G_B(\mathbb{R})$ of degree $d$ with multiple poles $\zeta_1, \cdots , \zeta_s$ with multiplicity $d_1, \cdots , d_s$, that is
\begin{equation*}
    \gamma = \frac{P(t)}{\prod_{i=1}^s (t + \zeta_i)^{d_i}(t - \zeta_i)^{d_i}}.
\end{equation*}
Theorem~\ref{decom_thm} states that $\gamma$ can be decomposed into a product of $s$ rational curves $\gamma_i$ of single pole $\zeta_i$ with multiplicity $d_i$, for $1 \le i \le s$. And the proof of Lemma~\ref{deg_red} in \cite{Liye2024rational} not only establishes this fact but also provides a constructive algorithm, which we present below in a form suitable for explicit computation. %
According to Lemma~\ref{deg_red}, there exist a rational curve
\begin{equation*}
    \alpha_1(t)= \frac{A(t)}{(t+ \zeta_1)^{d_1}(t - \zeta_1)^{d_1}},
\end{equation*}
where $A(t) = t^{2d_1} I_n + t^{2d_1-1} A_{2d_1 -1} + \cdots + A_0$, such that $\alpha_1\gamma$ is a rational curve of degree $d - 2d_1$.  Consequently, 
\begin{equation*}
    (t+ \zeta_1)^{2d_1}(t - \zeta_1)^{2d_1} \mid A(t)P(t),
\end{equation*}
which means $t = \zeta_1$ is a root of $A(t)P(t)$ of multiplicity $2d_1$, so the first $2d_1-1$ derivatives of $A(t)P(t)$ at point $t = \zeta_1$ are equal to $0$.  %
For $0 \leq l \leq 2d_1-1$, the $l$-th derivative of $A(t)P(t)$ at $t = \zeta_1$ is given by
 \begin{equation} \label{ini_lin}
    \left(AP\right)^{(l)}(\zeta_1) = \sum_{j = 0}^{l} \binom{l}{j} A^{(l-j)}(\zeta_1) P_j = 0.
 \end{equation}
 where $P_j = P^{(j)}(\zeta_1)$ for $0 \leq j \leq 2d_1-1$. 
And the $m$-th derivative of $A(t)$ at $t = \zeta_1$, for $0 \le m \le 2d_1-1$, is given by
\begin{equation*}
A^{(m)}(\zeta_1) = \frac{(2d_1)!}{(2d_1-m)!} \, I_n \, {\zeta_1}^{2d_1-m} + \sum_{j=m}^{2d_1-1} \frac{j!}{(j-m)!} A_j \, {\zeta_1}^{j-m}.
\end{equation*}
Thus we can rewrite \eqref{ini_lin} as the following linear system
\begin{equation} \label{mul_lin}
    \begin{pmatrix}
        I_n & A_{2d_1-1} & \cdots & A_1 & A_0
    \end{pmatrix} C_{\gamma} M_{\gamma}= 0.
\end{equation}
where 
\begin{equation*}
    M_{\gamma} = \begin{pmatrix}
        P_0 & \binom{1}{0} P_1 & \cdots & \binom{2d_1-1}{0} P_{2d_1-1} \\[3pt]
        0 & \binom{1}{1} P_0 &  \cdots &\binom{2d_1-1}{1} P_{2d_1-2} \\[3pt]
        0 & 0 &  \cdots & \binom{2d_1-1}{2} P_{2d_1-3} \\[3pt]
        \vdots  & \vdots & \ddots & \vdots \\[3pt]
        0 & 0 & \cdots & \binom{2d_1-1}{2d_1-1} P_0
    \end{pmatrix}
\end{equation*}
and
\begin{equation*}
    C_{\gamma} = \begin{pmatrix}
        \frac{(2d_1)! \, {\zeta_1}^{2d_1}}{(2d_1)!} I_n &  \cdots &  \frac{(2d_1)! \, {\zeta_1}^{1}}{(1)!} I_n \\[3pt]
        \frac{(2d_1-1)! \, {\zeta_1}^{2d_1-1}}{(2d_1-1)!} I_n  & \cdots &  \frac{(2d_1-1)! \, {\zeta_1}^{0}}{(0)!} I_n \\[3pt]
        \frac{(2d_1-2)! \, {\zeta_1}^{2d_1-2}}{(2d_1-2)!} I_n & \cdots &   0 \\[3pt]
        \vdots &  \ddots & \vdots \\[3pt]
        \frac{(0)! \, {\zeta_1}^{0}}{(0)!} I_n &  \cdots & 0
    \end{pmatrix}
\end{equation*}
Solving the linear system \eqref{mul_lin} yields a rational curve $\alpha_1$ with the property that rational curve $\alpha_1 \gamma$ has poles only at $\zeta_2, \dots, \zeta_s$ and its degree is reduced to $d - 2d_1$. Setting $\gamma_1 = \alpha_1^{-1}$, we obtain the first decomposition $\gamma = \gamma_1 (\alpha_1 \gamma)$.
We then apply the same procedure to $\alpha_1 \gamma$. This produces $\alpha_2$ such that rational curve $\alpha_2 \alpha_1 \gamma$ has poles only at $\zeta_3, \dots, \zeta_s$ and its degree becomes $d - 2d_1 - 2d_2$. Defining $\gamma_2 = \alpha_2^{-1}$, we refine the decomposition to $\gamma = \gamma_1 \gamma_2 (\alpha_2 \alpha_1 \gamma)$.
Proceeding iteratively in this manner, after $s$ steps we arrive at the decomposition
\begin{equation*}
   \gamma = \gamma_1 \cdots \gamma_s. 
\end{equation*}
where each $\gamma_i = \alpha_i^{-1}$ isolates the pole $\zeta_i$. By Corollary~\ref{uniq}, this decomposition is unique for a given order of the poles.
The complete procedure is formalized in Algorithm~\ref{alg:decomposition1}.
\begin{algorithm}[!htbp]
\normalsize\raggedright
\caption{Decompose multiple-pole rational curves}
\label{alg:decomposition1}
\setlist[enumerate]{leftmargin=*,topsep=3pt,itemsep=2pt}
\textbf{Input:} A rational curve on $G_B(\mathbb{R})$ 
\begin{equation*}
\gamma(t) = \frac{P(t)}{\prod_{i=1}^s (t + \zeta_i)^{d_i}(t - \zeta_i)^{d_i}},
\end{equation*}
where $P(t)$ is a polynomial matrix.
\par\smallskip\textbf{Output:} A decomposition into rational curves
\begin{equation*}
\gamma(t) = \gamma_1(t) \cdots \gamma_s(t) = \prod_{i=1}^s \frac{P_i(t)}{(t + \zeta_i)^{d_i}(t - \zeta_i)^{d_i}}.
\end{equation*}
\par\smallskip\hrule\smallskip
\textbf{For} $i = 1$ \textbf{to} $s$ \textbf{do}
\begin{enumerate}[label=\arabic*.,leftmargin=2.1em,labelsep=.45em]
    \item Define the intermediate term
    \[
    \alpha_i(t) = \frac{A_i(t)}{(t + \zeta_i)^{d_i}(t - \zeta_i)^{d_i}},
    \]
    where $A_i(t)$ is a polynomial matrix of degree $2d_i$:
    \[
    A_i(t) = t^{2d_i} I_n + t^{2d_i-1} A_{i,2d_i-1} + \dots + A_{i,0}.
    \]
    \item Determine $A_i(t)$ by solving the linear system:
    \begin{equation*}
    \left(
        I_n \quad A_{i,2d_i-1} \; \cdots \; A_{i,1} \quad A_{i,0}
   \right)
    C_{\gamma} M_{\gamma} = 0,
    \end{equation*}
    where $C_{\gamma}$ and $M_{\gamma}$ are matrices derived from $\gamma(t)$.
    \item Set $\gamma_i(t) = \alpha_i(t)^{-1}$.
    \item Update $\gamma(t) \leftarrow \alpha_i(t) \gamma(t)$.
\end{enumerate}
\textbf{End do}: \textbf{Return} $\gamma_1(t), \dots, \gamma_s(t)$.
\end{algorithm}
\begin{theorem}\label{alg1:correctness}
Algorithm~\ref{alg:decomposition1} terminates in finite steps 
and outputs.
\end{theorem}
\begin{proof}
Algorithm~\ref{alg:decomposition1} terminates in finite steps because the loop runs exactly $s$ times, where $s$ is the fixed number of distinct poles. The correctness of Algorithm~\ref{alg:decomposition1} follows from the analysis above and the proof of Lemma~\ref{deg_red} \cite{Liye2024rational}.
\end{proof}
The Table~\ref{tab:decomp-times}  shows that the decomposition time of rational curves on $SO_n(\mathbb{R})$ with different degree and different matrix dimension $n$. All decomposition times (total Maple kernel CPU time) were obtained using Maple 2022 on a laptop (Intel Core i5-8250U, 8 GB RAM).
\begin{table}[htbp]
\centering
\caption{Decomposition times of Algorithm~\ref{alg:decomposition1}}
\label{tab:decomp-times}
\begin{tabular}{lcc}
\toprule
 \textbf{Degree and Poles} & \textbf{Dimension} & \textbf{Decomposition Time (s)} \\
\midrule
d = 4  &$n=10$ & 0.610 \\
 poles = $\{\pm i, \pm 2 i\}$ &$n=20$ & 6.406 \\
 &$n=30$ & 31.391 \\
 &$n=50$ & 321.641 \\
 \midrule
d = 6 &$n=10$ & 1.422 \\
poles = $\{\pm i, \pm 2 i,  \pm 3 i\}$ &$n=20$ & 12.079 \\
 &$n=30$ & 60.859 \\
 &$n=50$ & 588.155 \\
\bottomrule
\end{tabular}
\end{table}

\FloatBarrier
\section{Decomposition Algorithm for Rational Curves with Single Pole}

\subsection{Decomposition with Extension}\label{subsec:ext}

While Algorithm~\ref{alg:decomposition1} provides a unique decomposition into single-pole curves, the resulting factors are not guaranteed to be of minimal degree. This leads to a natural subsequent question: can a rational curve with a single pole be further decomposed into a product of two lower-degree, single-pole curves? This problem is well-motivated by the factorization theory of motion polynomials \cite{GKLRSV17,Li25}, where a "multiplication trick" enables the reduction of degree when direct linear factorization is impossible. The feasibility of such a decomposition is supported by known examples. For instance,
\begin{equation*}
   \frac{\begin{pmatrix}
t^4 - 6t^2 + 1 & 4t^3 - 4t \\[4pt]
-4t^3 + 4t & t^4 - 6t^2 + 1
\end{pmatrix}}{(t^2+1)^2} = \frac{\begin{pmatrix}
t^{2} - 1 & 2t \\[4pt]
-2t & t^{2} - 1
\end{pmatrix}}{t^2+1}\frac{\begin{pmatrix}
t^{2} - 1 & 2t \\[4pt]
-2t & t^{2} - 1
\end{pmatrix}}{t^2+1}.
\end{equation*}
Analogous to our algorithm for decomposing rational curve with multiple poles, we can still reduce the degree of rational curves by solving the linear equations system. However, we can not guarantee that the solution is still rational curve, a counterexample as follows, which is given in \cite{Liye2024rational}
\begin{equation*}
\dfrac{
\begin{pmatrix}
2(t^2+1)^2 & 2t^2+2 & -2t^2-2 \\[4pt]
-2t^2-2 & 2t^4 + 4t^2 + 1 & 1 \\[4pt]
-2t^2-2 & -1 & 2t^4 + 4t^2 + 3
\end{pmatrix}}{2(t^2+1)^2}. 
\end{equation*}
Surprisingly, we observed that if we extended $\lambda$ onto $SO_{3,1}(\mathbb{R})$, 
then its can be decomposed into
\begin{equation*}
\dfrac{\begin{pmatrix}
8t^2 + 8 & 0 & 4t & -4t \\[4pt]
0 & 8t^2 + 8 & 4 & -4 \\[4pt]
-4t & -4 & 8t^2 + 7 & 1 \\[4pt]
-4t & -4 & -1 & 8t^2 + 9
\end{pmatrix}}{8(t^2+1)}\cdot \dfrac{ \begin{pmatrix}
8t^2 + 8 & 0 & -4t & 4t \\[4pt]
0 & 8t^2 + 8 & 4 & -4 \\[4pt]
4t & -4 & 8t^2 + 7 & 1 \\[4pt]
4t & -4 & -1 & 8t^2 + 9
\end{pmatrix}}{8(t^2+1)}. 
\end{equation*}
Therefore, in addition to solving the linear equation, we need to add quadratic equations to ensure that our solution is a rational curve. Nevertheless, we cannot guarantee that a solution exists for linear and quadratic systems. Thus, if no solution exists, we will extend the rational curve and then try to find a solution by solving the equations again. 
If $\lambda$ is a rational curve on $G_B(\mathbb{R})$ with single pole $\zeta$, without loss of generality, we can suppose $\zeta = i$ . Let $\lambda(t) = P(t)/(t^2+1)^d$, $d > 1$ be a rational curve on $SO_n$ or $O_{p,n-p}$, we say that rational curve $\lambda^e(t)$ is a extension of $\lambda(t)$ if $\lambda^e(t)$ is as follows
\begin{equation*}
    \lambda^e(t) = \dfrac{\begin{pmatrix}
        (t^2+1)^d & \mathbf{0} \\
        \mathbf{0} & P(t)
    \end{pmatrix}}{(t^2+1)^d}
    \quad \text{or} \quad
    \dfrac{\begin{pmatrix}
        P(t) & \mathbf{0} \\
        \mathbf{0} & \pm (t^2+1)^d 
    \end{pmatrix}}{(t^2+1)^d}.
\end{equation*}
\begin{remark}
    Although a rational curve can be extended to a higher order $G_B(\mathbb{R})$ in infinitely many ways, we focus solely on the canonical extensions presented here.
\end{remark}
Let $\alpha(t) = A(t)/(t^2+1)^s$, where $A(t) = t^{2s} I_n + \cdots + tA_1 + A_0$ and $s < d$. The condition that $(t^2+1)^{2s}$ divides $A(t)P(t)$ is equivalent to the linear system
\begin{equation} \label{2_lin}
    \begin{pmatrix}
        I_n & A_{2s-1} & \cdots & A_0  
    \end{pmatrix} C_{\lambda} M_{\lambda} = 0.
\end{equation}
Note that a solution $\alpha(t)$ of the linear system \eqref{2_lin} need not be a rational curve. To ensure that $\alpha(t)$ lies in the appropriate real algebraic group, we impose the additional quadratic system
\begin{equation} \label{quad_sys}
    A(t) B_{\lambda} A(t)^\mathsf{T} = (t^2+1)^{2s} B_{\lambda}.
\end{equation}
where $B_{\lambda}$ equals $I_n$ when $\lambda \in SO_n$, and $I_{p,n-p}$ when $\lambda \in O_{p,n-p}$.
If a solution $\alpha(t)$ to the systems \eqref{2_lin} and \eqref{quad_sys} exists, then $\alpha(t)$ has the desired properties: it is a rational curve on $SO_n$ or $O_{p,n-p}$, the product $\alpha\lambda$ has degree $2d-2s$, which is strictly lower than the original degree $2d$, and $\lambda$ admits the decomposition $\lambda = \alpha^{-1} \cdot (\alpha\lambda)$.
If no such solution exists for $\lambda$, we extend $\lambda$ to $\lambda^e$ and solve the extended linear system and quadratic system for $\lambda^e$ again. This extension introduces additional degrees of freedom, in hopes that the enlarged system becomes solvable.

\subsection{Necessary Condition}\label{sec:nece}
We establish a necessary condition for solvability based on the relationship between the ranks of the coefficient and augmented matrices for both the original and extended linear systems. 
%
%
Let $\lambda(t) = P(t)/(t^2+1)^d$, $d >1$ be a rational curve with single pole on $SO_n$ or $O_{p,n-p}$, we take $A(t) = t^{2s}I_n +  t^{2s-1} A_{2s-1} + \cdots + A_0$ for $ s < d$,  where $A_{2s-1}, \cdots, A_0$ are $n \times n$ real variables matrices, then consider the following linear system
\begin{equation}\label{linear_eq}
\begin{pmatrix}
        I_n & A_{2s-1} & \cdots & A_0  
    \end{pmatrix} C_{\lambda} M_{\lambda} = 0.
\end{equation}
where \begin{equation*}
    M_{\lambda} = \begin{pmatrix}
         P_0 & \binom{1}{0} P_1  & \cdots & \binom{2s-1}{0} P_{2s-1} \\[3pt]
         0 & \binom{1}{1} P_0  & \cdots &\binom{2s-1}{1} P_{2s-2} \\[3pt]
        \vdots &\vdots  & \ddots & \vdots \\[3pt]
         0 & 0 &   \cdots & \binom{2s-1}{2s-1} P_0
     \end{pmatrix}
 \end{equation*}
 here $P_j = P^{(j)}(i)$ for $0 \leq j \leq 2s-1$, and
 \begin{equation*}
     C_{\lambda} = \begin{pmatrix}
         \frac{(2s)! \, i^{2s}}{(2s)!} I_n & \cdots &  \frac{(2s)! \, i^{1}}{(1)!} I_n \\[3pt]
         \frac{(2s-1)! \, i^{2s-1}}{(2s-1)!} I_n  & \cdots &  \frac{(2s-1)! \, i^{0}}{(0)!} I_n \\[3pt]
        
         \vdots &  \ddots & \vdots \\[3pt]
         \frac{(0)! \, i^{0}}{(0)!} I_n &   \cdots & 0
    \end{pmatrix}
 \end{equation*}
Since we'll find real solutions to the complex linear system \eqref{linear_eq}, we first reformulate it as an equivalent real linear system. Thus the solvability of \eqref{linear_eq} then depends on the rank properties of the corresponding real coefficient and augmented matrices.
Let $S_{\lambda} := (C_{\lambda}M_{\lambda})^\mathsf{T}$, $S_{\lambda}$ is $2sn \times (2s+1)n$ matrix. Take the $(n+1)$ to $(2s+1)n$ columns of $S_{\lambda}$ to form a $2sn \times 2sn$ matrix, denoted by $K_{\lambda}$, and let $L_{\lambda}$ be a diagonal matrix composed of n matrices $K_{\lambda}$
\begin{equation*}
L_{\lambda} = \begin{pmatrix}
K_{\lambda}  & & \\
&  \ddots & \\
& & K_{\lambda}
\end{pmatrix}_{2sn^2 \times 2sn^2}
\end{equation*}
Then, we separate  $L_{\lambda}$ into real and imaginary parts. Let $L_{\mathrm{re}}$ and $L_{\mathrm{im}}$ be the real and imaginary parts of $L_{\lambda}$ respectively. These are combined into the real coefficient matrix
\begin{equation*}
T_{\lambda} =\begin{pmatrix}
L_{\mathrm{re}}\\
L_{\mathrm{im}}
\end{pmatrix}_{4sn^2 \times 2sn^2}
\end{equation*}
For the constant term in the real linear system, take the first $n$ columns of $S_{\lambda}=(s_{ij})$, where $ 1\le i \le 2sn$, $1 \le j \le (2s+1)n$, to form the complex vector in $ \mathbb{R}^{2sn^2}$
\begin{equation*}
v_{\lambda} = \left(
s_{11}  \cdots  s_{2sn,1} \; s_{12}  \cdots s_{2sn,2}\cdots \cdots s_{1,n} \cdots s_{2sn,n}
\right)^\mathsf{T}.
\end{equation*}
Let $v_{\mathrm{re}}$ and $v_{\mathrm{im}}$ be the real and imaginary parts of $v_{\lambda}$. And the real vector $b_{\lambda} \in \mathbb{R}^{4sn^2}$ is then defined as their concatenation:
\begin{equation*}
b_{\lambda} = \begin{pmatrix}
    v_{\mathrm{re}}\\
    v_{\mathrm{im}}
\end{pmatrix}.
\end{equation*}
Consequently, the complex linear system \eqref{linear_eq} is equivalent to the real linear system 
\begin{equation*}
   T_{\lambda} X = -b_{\lambda}. 
\end{equation*}
where $X \in \mathbb{R}^{2sn^2}$ consists of the entries of the matrices $A_{2s-1}, \cdots, A_0$. Explicitly, if $A_k = (a^k_{ij})_{n \times n}$ for $1 \le k \le 2s-1$, then 
 \begin{equation*}
  X = \left(
a^{2s-1}_{11} \cdots a^{2s-1}_{1n} \cdots  a^0_{11}  \cdots a^0_{1n}\cdots  a^{2s-1}_{n1} \cdots a^{2s-1}_{nn} \cdots a^0_{n1}  \cdots a^0_{nn}
\right)^\mathsf{T}.
 \end{equation*}
 \begin{remark}
     Note that the matrices and vectors introduced above depend on both the rational curve $\lambda$ and the integer parameter $s$. To keep notation light, we suppress the dependence on $s$, writing $T_{\lambda}$, $b_{\lambda}$, etc., in place of the more precise $T_{\lambda,s}$, $b_{\lambda,s}$, unless explicit clarification is needed.
 \end{remark}
\begin{lemma} \label{eq}
Let $\lambda(t)=P(t)/(t^2+1)^d$ with $d>1$ be a rational curve. For integer $s < d$, there exists a matrix polynomial 
\begin{equation*}
  A(t) = t^{2s}I_n +  t^{2s-1} A_{2s-1} + \cdots + A_0, \; A_k \in \mathbb{R}^{n \times n},
\end{equation*}
such that $(t^2+1)^{2s}$ divides $A(t)P(t)$ if and only if
\begin{equation*}
  \operatorname{rank}(T_{\lambda})=\operatorname{rank}\bigl(T_{\lambda}\mid -b_{\lambda}\bigr).
\end{equation*}
\end{lemma}
\begin{proof}
By construction, the existence of a matrix $A(t)$ with the stated property is equivalent to the solvability of the linear system $T_{\lambda}X = -b_{\lambda}$. And this real linear system has a solution if and only if the rank of its coefficient matrix equals the rank of augmented matrix, that is, $\operatorname{rank}(T_{\lambda})=\operatorname{rank}\bigl(T_{\lambda}\mid -b_{\lambda}\bigr)$. The equivalence stated in the lemma follows immediately.
\end{proof}
To prepare for the proof of Theorem~\ref{nece}, we first establish a useful relation between the ranks of $T_\lambda$ and its block component. Let $K_{\mathrm{re}}$ and $K_{\mathrm{im}}$ denote the real and imaginary parts of $K_{\lambda}$, and define
\begin{equation*}
Q_{\lambda} = \begin{pmatrix}
K_{\mathrm{re}} \\
K_{\mathrm{im}}
\end{pmatrix}_{4sn \times 2sn}
\end{equation*}

\begin{lemma}\label{nrank}
The ranks satisfy $\operatorname{rank}(T_{\lambda}) = n \cdot \operatorname{rank}(Q_{\lambda})$.
\end{lemma}
\begin{proof}
By performing elementary row operations that do not alter the rank, we can transform $T_{\lambda}$ into  a block-diagonal form
\begin{equation*}
\begin{pmatrix}
     Q_{\lambda}  & & \\
&  \ddots & \\
&  & Q_{\lambda}
\end{pmatrix}_{4sn^2 \times 2sn^2}
\end{equation*}
Since the rank of a block diagonal matrix is the sum of the ranks of its diagonal blocks, we obtain $\operatorname{rank}(T_{\lambda}) = n \cdot \operatorname{rank}(Q_{\lambda})$.
\end{proof}
%
%
%
%
\begin{theorem}\label{nece}
Let $\lambda(t) = P(t)/(t^2+1)^d$, $d>1$ be a rational curve, and $\lambda^e$ be an extension of $\lambda$. For integer $s$ satisfying $2s-1 < d$, then
\begin{equation*}
  \operatorname{rank}\bigl(T_{\lambda} \mid -b_{\lambda}\bigr) - \operatorname{rank}(T_{\lambda})
= \operatorname{rank}\bigl(T_{\lambda^e} \mid -b_{\lambda^e}\bigr) - \operatorname{rank}(T_{\lambda^e}).  
\end{equation*}
\end{theorem}
\begin{proof}
Without loss of generality, we take the extension to be
\begin{equation*}
\lambda^e(t) = \frac{P^e(t)}{(t^2+1)^d}, \; \text{where} \; P^e(t) = \begin{pmatrix}
(t^2+1)^d & 0 \\
0 & P(t)
\end{pmatrix}. 
\end{equation*}
Since the derivatives of $(t^2+1)^d$ at $t=i$ vanish up to order $d-1$, we have
\begin{equation*}
P^e_j := \bigl(P^e\bigr)^{(j)}(i) = \begin{pmatrix}
0 & 0 \\
0 & P_j
\end{pmatrix}.
\end{equation*}
where $P_j = P^{(j)}(i)$ for $ 0 \le j \le 2s-1$. 
Consequently, $S_{\lambda^e}$ is obtained from $S_{\lambda}$ by adjoining zero rows and columns; the same holds for $K_{\lambda^e}$ and $Q_{\lambda^e}$. These modifications do not affect the rank, hence $\operatorname{rank}(Q_{\lambda^e}) = \operatorname{rank}(Q_{\lambda})$. And by Lemma~\ref{nrank}, we have 
\begin{equation}\label{rela}
\begin{aligned} 
\operatorname{rank}(T_{\lambda^e}) & = (n+1) \cdot \operatorname{rank}(Q_{\lambda^e}) \\
& = (n+1)\cdot  \operatorname{rank}(Q_{\lambda}) \\
& =  \operatorname{rank}(Q_{\lambda}) + \operatorname{rank}(T_{\lambda}).
\end{aligned}  
\end{equation}
For the constant term, a similar inspection shows
\begin{equation*}
v_{\lambda^e} =\begin{pmatrix}
\mathbf{0} \\
v'_{\lambda}
\end{pmatrix}.
\end{equation*}
where $v'_{\lambda}$ is obtained from $v_{\lambda}$ by inserting appropriate zeros, and $\mathbf{0}$ is $2sn$-ary zero vector.
Performing elementary row operations that preserve the rank, we transform augmented matrix $(T_{\lambda^e} \mid -b_{\lambda^e})$ into
\begin{equation*}
   \begin{pmatrix}
Q_{\lambda^e} &        &        & \mathbf{0} \\
              & \ddots &        & \vdots \\
              &        & Q_{\lambda^e} & b'_{\lambda}
\end{pmatrix}
\end{equation*}
where $b'_{\lambda}$ is $b_{\lambda}$ with some zeros added, and $\mathbf{0}$ is $4sn$-ary zero vector.
Thus we have 
\begin{equation} \label{rela'}
\operatorname{rank}\bigl(T_{\lambda^e} \mid -b_{\lambda^e}\bigr) 
= \operatorname{rank}(Q_{\lambda^e}) + \operatorname{rank}\bigl(T_{\lambda} \mid -b_{\lambda}\bigr).
\end{equation}
and by \eqref{rela} and \eqref{rela'}, we have
\begin{equation*}
   \operatorname{rank}\bigl(T_{\lambda} \mid -b_{\lambda}\bigr) - \operatorname{rank}(T_{\lambda})
= \operatorname{rank}\bigl(T_{\lambda^e} \mid -b_{\lambda^e}\bigr) - \operatorname{rank}(T_{\lambda^e}).  
\end{equation*}
\end{proof}

\subsection{Decomposition Algorithm}\label{sec:decomp2}
Based on the necessary condition derived above for the linear system, we present an algorithm for decomposing rational curves with single pole.
Let $\lambda = P(t)/(t^2+1)^d$, $d > 1$ be a rational curve on $SO_n$ or $O_{p,n-p}$, given a candidate parameter $s$ satisfying $2s - 1 < d$, and let $\alpha = A(t)/(t^2+1)^s$,  the procedure is as follows.
First, we compute the ranks of the coefficient matrix $T_{\lambda,s}$ and the augmented matrix  $(T_{\lambda,s}\mid b_{\lambda,s})$. If these ranks are equal, then we attempt to solve the linear and quadratic systems to obtain a  decomposition.
If the ranks differ, it indicates that no solution exists for the current choice of $s$. In this case, increment the candidate $s$ and re-evaluates the ranks of linear system with the updated parameter.
If the ranks are equal but the linear and quadratic system admits no solution, this suggests that the target decomposition cannot be achieved for the original rational curve. Then we extend the rational curve $\lambda$ to  $\lambda^e$ and solve the extended system iteratively .
Finally, if no decomposition exists for any admissible $s$ and  any extended rational curve $\lambda^e$ , then output Null. This signifies that the given rational curve $\lambda$ is indecomposable by the present algorithm.
 %
 %
\begin{remark}
    The condition $2s-1 < d$ is required since Theorem~\ref{nece} generally holds for $2s-1 < d$.
\end{remark}
The complete procedure is formalized in Algorithm~\ref{alg:decomposition2}.
\begin{algorithm}[!htbp]
\normalsize\raggedright
\caption{Decompose single-pole rational curves}
\label{alg:decomposition2}
\setlist{leftmargin=*,topsep=3pt,itemsep=2pt}
\textbf{Input:} A rational curve $\lambda = P(t)/(t^2+1)^d$ on $SO_n \;\text{or} \; O_{p,n-p}$
where $d > 1$ and $P(t)$ is a polynomial matrix.
%
%
%
\par\smallskip\textbf{Output:}
\begin{itemize}
    \item \textbf{Success case:} A decomposition of $\lambda$ into: $\lambda = \lambda_1 \lambda_2 = (P_1(t)/(t^2+1)^s) \cdot (P_2(t)/(t^2+1)^{d-s})$

or $\lambda^e = \lambda_1 \lambda_2$ in the extended case.
    \item \textbf{Failure case:} \texttt{Null}.
\end{itemize}

\par\smallskip\hrule\smallskip

\textbf{Initialize} $s \leftarrow 1$.

\textbf{While} $2s-1 < d$ \textbf{do}\par

\begin{algorithmindent}

 \textbf{If} $\operatorname{rank}(T_{\lambda,s})=\operatorname{rank}\bigl(T_{\lambda,s}\mid -b_{\lambda,s}\bigr)$ \textbf{then}

\begin{enumerate}[label=\arabic*.,leftmargin=2.1em,labelsep=.45em]
    \item \begin{minipage}[t]{\linewidth}\raggedright
    \textbf{Attempt decomposition in the original space.}

    Let $\alpha(t)=A(t)/(t^2+1)^s$, $A(t)=t^{2s}I_n+\sum_{j=0}^{2s-1}t^j A_j$.

    Solve the following system for $A_{2s-1}, \cdots, A_0$:
    \[
    \begin{cases}
        \left( I_n \; A_{2s-1} \; \cdots \; A_1  \; A_0 \right) C_{\lambda,s} M_{\lambda,s} = 0, \\[4pt]
        A(t) B_{\lambda} A(t)^\mathsf{T} = (t^2 + 1)^{2s} B_{\lambda}.
    \end{cases}
    \]

    \textbf{If} solvable \textbf{then} set $\lambda_1 = \alpha^{-1}$, $\lambda_2 = \alpha \lambda$, \textbf{return} $\lambda_1, \lambda_2$.

\textbf{Else}  proceed to Step~2.\end{minipage}

    \item \begin{minipage}[t]{\linewidth}\raggedright
    \textbf{Attempt decomposition in the extended space.}

    Set $\lambda(t) \leftarrow \lambda^e(t)$.

    Let $\alpha(t)=A(t)/(t^2+1)^s$, $A(t)=t^{2s}I_{n+1}+\sum_{j=0}^{2s-1}t^j A_j$.

    Solve the corresponding system of $\lambda^e$.

    \textbf{If} solvable \textbf{then} set $\lambda_1 = \alpha^{-1}$, $\lambda_2 = \alpha \lambda^e$, \textbf{return} $\lambda_1, \lambda_2$.

    \textbf{Else} set $s \leftarrow s + 1$.\end{minipage}
\end{enumerate}

\noindent\textbf{Else} (rank condition fails): set $s \leftarrow s + 1$.

\end{algorithmindent}


    
    
    

   
    
    

    



\textbf{End do}: \textbf{Return} \texttt{Null}. 

\end{algorithm}

\begin{theorem}\label{alg2:correctness}
Algorithm~\ref{alg:decomposition2} terminates in finite steps 
and outputs.
\end{theorem}

\begin{proof}
Algorithm~\ref{alg:decomposition2} terminates in finite steps because the loop condition is $2s - 1 < d$, thus after at most 
$\lfloor (d+1)/2 \rfloor$ iterations, the loop ends, either by returning a decomposition or Null. The correctness of Algorithm~\ref{alg:decomposition2} follows from Theorem~\ref{nece} and Lemma~\ref{deg_red}.
\end{proof}

\section{Examples}\label{sec:examples}

We now demonstrate our algorithms on key examples. Example~\ref{example1} illustrates Algorithm~\ref{alg:decomposition1}. Examples~\ref{ex:translation} and ~\ref{example3} demonstrate Algorithm~\ref{alg:decomposition2}, showing successful decomposition and the role of the necessary rank condition.

\begin{example}\label{example1}
Let $\lambda$ be a rational curve on $SO_3(\mathbb{R})$ of degree $4$ and poles $\{\pm i, \pm 2i\}$ as follows
    \begin{equation*}
\dfrac{\begin{pmatrix}
t^{4} - 3t^{2} - 4 & 4t^{3} + 4t & 0 \\[3pt]
-4t^{3} + 4t & t^{4} - 5t^{2} + 4 & 2t^{3} + 8t \\[3pt]
8t^{2} & -2t^{3} + 8t & t^{4} + 3t^{2} - 4
\end{pmatrix}}{(t^2+1)(t^2+4)}.
\end{equation*}
By Algorithm~\ref{alg:decomposition1}, $\lambda$ can be decomposed into
    \begin{equation*}
\dfrac{\begin{pmatrix}
t^2 + 1 & 0 & 0 \\
0 & t^2 - 1 & 2t \\
0 & -2t & t^2 - 1
\end{pmatrix}}{(t^2+1)}\cdot\dfrac{\begin{pmatrix}
t^2 - 4 & 4t & 0 \\
-4t & t^2 - 4 & 0 \\
0 & 0 & t^2 + 4
\end{pmatrix}}{(t^2+4)}
\end{equation*}
or
 \begin{equation*}
\dfrac{\begin{pmatrix}
25t^{2} + 28 & 60t & 96 \\[4pt]
-60t & 25t^{2} - 100 & 80t \\[4pt]
96 & -80t & 25t^{2} - 28
\end{pmatrix}}{25(t^2+4)}\cdot\dfrac{\begin{pmatrix}
25t^{2} - 7 & 40t & -24 \\[4pt]
-40t & 25t^{2} - 25 & -30t \\[4pt]
-24 & 30t & 25t^{2} + 7
\end{pmatrix}}{25(t^2+1)}
\end{equation*}
\end{example}

\begin{example}\label{ex:translation}
Let $\lambda$ be a rational curve on $\text{SO}_{4,1}$ of degree $4$ as follows
 \begin{equation*}
    \dfrac{\begin{pmatrix}
2(t^2 + 1)^2 & 0 & 0 & 0 & 0 \\[4pt]
0 & 2(t^2 + 1)^2 & 0 & 2(t^3 + t) & -2(t^3 + t) \\[4pt]
0 & 0 & 2(t^2 + 1)^2 & 0 & 0 \\[4pt]
0 & -2(t^3 + t) & 0 & 2(t^2+1)^2 - t^2 & t^2 \\[4pt]
0 & -2(t^3 + t) & 0 & -t^2 & 2(t^2+1)^2 + t^2
\end{pmatrix}}{2(t^2+1)^2}
    \end{equation*}
    Then $\lambda$ is obtained using  a translation motion polynomial $t^2+1+\epsilon \mathrm{i}$ defined in \cite{Hegedus13, Thimm2026}, which can be mapped to a quadratic curve on the special Euclidean group SE$(3)$ as
    \begin{equation*}
    \begin{pmatrix}
1 & 0 & 0 & 0 \\
0 & 1 & 0 & \frac{2t}{t^2+1} \\
0 & 0 & 1 & 0 \\
0 & 0 & 0 &  1
\end{pmatrix}
    \end{equation*}
    Using Algorithm~\ref{alg:decomposition2},  $\lambda$ can be decomposed into the product of the two quadratic curves, and the decomposition is not unique 


\begin{equation*}
   \dfrac{\begin{pmatrix}
2K(t^2+1) & 0 & 0 & 8u(tP+vQ) & -8u(tP+vQ) \\[4pt]
0 & 2K(t^2+1) & 0 & Kt & -Kt \\[4pt]
0 & 0 & 2K(t^2+1) & 8u(Q-tvP) & 8u(tvP-Q) \\[4pt]
-8u(tP+vQ) & -Kt & 8u(tvP-Q) & 2K(t^2+1)-Q^2 & Q^2 \\[4pt]
-8u(tP+vQ) & -Kt & 8u(tvP-Q) & -Q^2 & 2K(t^2+1)+Q^2
\end{pmatrix}}{2K (t^2+1)}
\end{equation*}
\begin{equation*}
   \dfrac{\begin{pmatrix}
2K(t^2+1) & 0 & 0 & -8u(tP+vQ) & 8u(tP+vQ) \\[4pt]
0 & 2K(t^2+1) & 0 & Kt & -Kt \\[4pt]
0 & 0 & 2K(t^2+1) & -8u(Q-tvP) & -8u(tvP-Q) \\[4pt]
8u(tP+vQ) & -Kt & -8u(tvP-Q) & 2K(t^2+1)-Q^2 & Q^2 \\[4pt]
8u(tP+vQ) & -Kt & -8u(tvP-Q) & -Q^2 & 2K(t^2+1)+Q^2
\end{pmatrix}}{2K (t^2+1)}
\end{equation*}
where $P = u^2 - v^2 - 1, \, Q = u^2 + v^2 + 1, \, R = u(v^2+1), \, K = 16 u^2(v^2+1)$.
\end{example}

\begin{example}\label{example3}

Let $\lambda$ be a rational curve on $\text{SO}_{4,1}$ of degree $8$ as follows

\begin{equation*}
    \dfrac{\begin{pmatrix}
2(t^2 + 1)^4 & 0 & 0 & 0 & 0 \\[4pt]
0 & 2(t^2 + 1)^4 & 0 & 2(t^2 + 1)^2 t & -2(t^2 + 1)^2 t \\[4pt]
0 & 0 & 2(t^2 + 1)^4 & 0 & 0 \\[4pt]
0 & -2(t^2 + 1)^2 t & 0 & 2(t^2+1)^4 - t^2 & t^2 \\[4pt]
0 & -2(t^2 + 1)^2 t & 0 & -t^2 & 2(t^2+1)^4 + t^2
\end{pmatrix}}{2(t^2+1)^4}
\end{equation*}
Compute the rank of $\, T_{\lambda,s}$ and $(T_{\lambda,s} \mid -b_{\lambda,s})$ for $s=1, 2$, we have
\begin{equation*}
\text{rank}\,(T_{\lambda,1}) = 10 \neq \text{rank}\,(T_{\lambda,1} \mid -b_{\lambda,1}) = 11.
\end{equation*}
and
\begin{equation*}
\text{rank}\,(T_{\lambda,2}) =  \text{rank}\,(T_{\lambda,2} \mid -b_{\lambda,2}) = 30.
\end{equation*}
Then by Algorithm \ref{alg:decomposition2}, $\lambda$ can be decomposed into the product of the following lower degree rational curves
\begin{equation*}
    \dfrac{\begin{pmatrix}
    16(t^2 + 1)^2 & 0 & 0 & 4\sqrt{2} & -4\sqrt{2} \\[4pt]
    0 & 16(t^2 + 1)^2 & 0 & 8t & -8t \\[4pt]
    0 & 0 & 16(t^2 + 1)^2 & 4\sqrt{2}t^{2} & -4\sqrt{2}t^{2} \\[4pt]
    -4\sqrt{2} & -8t & -4\sqrt{2}t^{2} & 16(t^2+1)^2 - 1 & 1 \\[4pt]
    -4\sqrt{2} & -8t & -4\sqrt{2}t^{2} & -1 & 16(t^2+1)^2 + 1
    \end{pmatrix}}{16(t^2+1)^2}
\end{equation*}


\begin{equation*}
    \dfrac{\begin{pmatrix}
16(t^2 + 1)^2 & 0 & 0 & -4\sqrt{2} & 4\sqrt{2} \\[4pt]
0 & 16(t^2 + 1)^2 & 0 & 8t & -8t \\[4pt]
0 & 0 & 16(t^2 + 1)^2 & -4\sqrt{2}t^{2} & 4\sqrt{2}t^{2} \\[4pt]
4\sqrt{2} & -8t & 4\sqrt{2}t^{2} & 16(t^2+1)^2 - 1 & 1 \\[4pt]
4\sqrt{2} & -8t & 4\sqrt{2}t^{2} & -1 & 16(t^2+1)^2 + 1
\end{pmatrix}}{16(t^2+1)^2}
\end{equation*}

\end{example}

\section{Summary}\label{sec:summary}
This paper presents two novel decomposition algorithms for rational curves on real classical matrix groups, 
$G_B(\mathbb{F})$, with complete implementations in Maple. The Algorithm~\ref{alg:decomposition1} uniquely decomposes a rational curve with multiple poles into a product of curves each having only a single pole, solving only linear equations. The Algorithm~\ref{alg:decomposition2}, refined algorithm, tackles the more complex problem of decomposing a single pole curve into lower degree factors by solving a system of linear and quadratic equations, a process motivated by and shown to generalize the "multiplication trick" in motion polynomial factorization.
A key theoretical contribution is the proof that the solvability of the underlying linear system is a necessary condition for Algorithm~\ref{alg:decomposition2} to succeed; if this condition fails, no solution exists even after extending the system. While a sufficient condition for guaranteed decomposition remains an open theoretical question, the algorithm demonstrated empirical success on all tested examples. The work provides practical tools for symbolic computation in geometric algebra and clarifies the relationship between curve decomposition in matrix models and classical factorization problems in kinematics.





	\bibliographystyle{plainnat}
	\bibliography{rationalcurves}

@article{Li25,
  author    = {Zijia Li and Hans-Peter Schr{\"o}cker and Johannes Siegele},
  title     = {A geometric algorithm for the factorization of rational motions in conformal three space},
  journal   = {Journal of Symbolic Computation},
  year      = {2025},
  volume    = {128},
  pages     = {102388},
  doi       = {10.1016/j.jsc.2024.102388},
  issn      = {0747-7171}
}

@Article{Hegedus13,
  author =       {Heged\"us, G\'abor and Schicho, Josef and Schr\"ocker, Hans-Peter},
  title =        {Factorization of Rational Curves in the {Study} Quadric},
  journal =      {Mech. Machine Theory},
  year =         2013,
  volume =       69,
  number =       1,
  pages =        {142--152},
  doi =    {10.1016/j.mechmachtheory.2013.05.010},
}

@article{Liye2024rational,
  title={Rational Curves on Real Classical Groups},
  author={Zijia Li and Ke Ye},
  journal={ArXiv},
  year={2024},
  volume={abs/2408.04453},
  url={https://api.semanticscholar.org/CorpusID:271768708}
}

@article{Juttler93zwanglaufige,
  title={{\"U}ber zwangl{\"a}ufige rationale Bewegungsvorg{\"a}nge},
  author={J{\"u}ttler, Bert},
  journal={{\"O}sterreich. Akad. Wiss. Math.-Natur. Kl. S.-B. II},
  volume={202},
  number={1-10},
  pages={117--232},
  year={1993}
}

@book {Winkler08,
    AUTHOR = {Sendra, J. Rafael and Winkler, Franz and P\'erez-D\'iaz,
              Sonia},
     TITLE = {Rational algebraic curves},
    SERIES = {Algorithms and Computation in Mathematics},
    VOLUME = {22},
      NOTE = {A computer algebra approach},
 PUBLISHER = {Springer, Berlin},
      YEAR = {2008},
     PAGES = {x+267},
      ISBN = {978-3-540-73724-7},
   MRCLASS = {14Q05 (14H50 68W30)},
  MRNUMBER = {2361646},
MRREVIEWER = {Valmecir\ A.\ Bayer},
       DOI = {10.1007/978-3-540-73725-4},
       URL = {https://doi.org/10.1007/978-3-540-73725-4},
}

@book {Roth79tk,
    AUTHOR = {Bottema, O. and Roth, B.},
     TITLE = {Theoretical kinematics},
    SERIES = {North-Holland Series in Applied Mathematics and Mechanics},
    VOLUME = {24},
 PUBLISHER = {North-Holland Publishing Co., Amsterdam-New York},
      YEAR = {1979},
     PAGES = {xiv+558},
      ISBN = {0-444-85124-0},
   MRCLASS = {70-02 (51Nxx 53A17 70B99)},
  MRNUMBER = {533960},
MRREVIEWER = {Robert\ Connelly},
}

@article {Darboux1890,
    AUTHOR = {Darboux, Gaston},
     TITLE = {Sur le d\'eplacement d'une figure invariable},
   JOURNAL = {Ann. Sci. \'Ecole Norm. Sup. (3)},
  FJOURNAL = {Annales Scientifiques de l'\'Ecole Normale Sup\'erieure.
              Troisi\`eme S\'erie},
    VOLUME = {7},
      YEAR = {1890},
     PAGES = {323--326},
      ISSN = {0012-9593},
   MRCLASS = {99-04},
  MRNUMBER = {1508846},
       URL = {http://www.numdam.org/item?id=ASENS_1890_3_7__323_0},
}

@article {Kollar92fano,
    AUTHOR = {Koll\'ar, J\'anos and Miyaoka, Yoichi and Mori, Shigefumi},
     TITLE = {Rational connectedness and boundedness of {F}ano manifolds},
   JOURNAL = {J. Differential Geom.},
  FJOURNAL = {Journal of Differential Geometry},
    VOLUME = {36},
      YEAR = {1992},
    NUMBER = {3},
     PAGES = {765--779},
      ISSN = {0022-040X,1945-743X},
   MRCLASS = {14J45},
  MRNUMBER = {1189503},
MRREVIEWER = {Yuri\ G.\ Prokhorov},
       URL = {http://projecteuclid.org/euclid.jdg/1214453188},
}

@article {Kollar92,
    AUTHOR = {Koll\'ar, J\'anos and Miyaoka, Yoichi and Mori, Shigefumi},
     TITLE = {Rationally connected varieties},
   JOURNAL = {J. Algebraic Geom.},
  FJOURNAL = {Journal of Algebraic Geometry},
    VOLUME = {1},
      YEAR = {1992},
    NUMBER = {3},
     PAGES = {429--448},
      ISSN = {1056-3911,1534-7486},
   MRCLASS = {14E30 (14E05 14E35)},
  MRNUMBER = {1158625},
MRREVIEWER = {Yuri\ G.\ Prokhorov},
}

@Article{Mori1988flip,
  title={Flip theorem and the existence of minimal models for 3-folds},
  author={Mori, Shigefumi},
  journal={Journal of the American Mathematical Society},
  volume={1},
  number={1},
  pages={117--253},
  year={1988},
  publisher={JSTOR}
}

@article {Graber03families,
    AUTHOR = {Graber, Tom and Harris, Joe and Starr, Jason},
     TITLE = {Families of rationally connected varieties},
   JOURNAL = {J. Amer. Math. Soc.},
  FJOURNAL = {Journal of the American Mathematical Society},
    VOLUME = {16},
      YEAR = {2003},
    NUMBER = {1},
     PAGES = {57--67},
      ISSN = {0894-0347,1088-6834},
   MRCLASS = {14M20 (14D05 14D22)},
  MRNUMBER = {1937199},
MRREVIEWER = {Alexandr\ V.\ Pukhlikov},
       DOI = {10.1090/S0894-0347-02-00402-2},
       URL = {https://doi.org/10.1090/S0894-0347-02-00402-2},
}

@article{Thimm2026,
  author    = {Daren A. Thimm and Zijia Li and Hans-Peter Schr{\"o}cker and Johannes Siegele},
  title     = {Quadratic Motion Polynomials with Irregular Factorizations},
  journal   = {Advances in Applied Clifford Algebras},
  year      = {2026},
  volume    = {36},
  number    = {1},
  pages     = {19},
  doi       = {10.1007/s00006-025-01426-2},
  issn      = {1661-4909},
  date      = {2026-01-24}
}

@article {LJS18,
    AUTHOR = {Li, Zijia and Schicho, Josef and Schr\"{o}cker, Hans-Peter},
     TITLE = {Kempe's universality theorem for rational space curves},
   JOURNAL = {Found. Comput. Math.},
  FJOURNAL = {Foundations of Computational Mathematics. The Journal of the
              Society for the Foundations of Computational Mathematics},
    VOLUME = {18},
      YEAR = {2018},
    NUMBER = {2},
     PAGES = {509--536},
      ISSN = {1615-3375,1615-3383},
   MRCLASS = {70B05 (14H50 65D17 68U07)},
  MRNUMBER = {3777787},
MRREVIEWER = {Peter\ S.\ Donelan},
       DOI = {10.1007/s10208-017-9348-x},
       URL = {https://doi.org/10.1007/s10208-017-9348-x},
}

@article {KM02,
    AUTHOR = {Kapovich, Michael and Millson, John J.},
     TITLE = {Universality theorems for configuration spaces of planar
              linkages},
   JOURNAL = {Topology},
  FJOURNAL = {Topology. An International Journal of Mathematics},
    VOLUME = {41},
      YEAR = {2002},
    NUMBER = {6},
     PAGES = {1051--1107},
      ISSN = {0040-9383},
   MRCLASS = {14P05 (55R80)},
  MRNUMBER = {1923214},
MRREVIEWER = {Thomas\ C.\ Craven},
       DOI = {10.1016/S0040-9383(01)00034-9},
       URL = {https://doi.org/10.1016/S0040-9383(01)00034-9},
}

@article {KM98,
    AUTHOR = {Kapovich, Michael and Millson, John J.},
     TITLE = {On representation varieties of {A}rtin groups, projective
              arrangements and the fundamental groups of smooth complex
              algebraic varieties},
   JOURNAL = {Inst. Hautes \'{E}tudes Sci. Publ. Math.},
  FJOURNAL = {Institut des Hautes \'{E}tudes Scientifiques. Publications
              Math\'{e}matiques},
    Volume = {88},
      YEAR = {1998},
     PAGES = {5--95},
      ISSN = {0073-8301,1618-1913},
   MRCLASS = {14F35 (14B12 14C30 14D20 14E20 20F36)},
  MRNUMBER = {1733326},
MRREVIEWER = {Philippe\ P.\ Eyssidieux},
       URL = {http://www.numdam.org/item?id=PMIHES_1998__88__5_0},
}

@article {KR15,
    AUTHOR = {Kharlamov, Viatcheslav and R\u{a}sdeaconu, Rare\c{s}},
     TITLE = {Counting real rational curves on {$K3$} surfaces},
   JOURNAL = {Int. Math. Res. Not. IMRN},
  FJOURNAL = {International Mathematics Research Notices. IMRN},
      YEAR = {2015},
    NUMBER = {14},
     PAGES = {5436--5455},
      ISSN = {1073-7928,1687-0247},
   MRCLASS = {14N10 (14J28 14P05)},
  MRNUMBER = {3384445},
MRREVIEWER = {Hassan\ Haghighi},
       DOI = {10.1093/imrn/rnu091},
       URL = {https://doi.org/10.1093/imrn/rnu091},
}

@article {KR17,
    AUTHOR = {Kharlamov, Viatcheslav and R\u{a}sdeaconu, Rare\c{s}},
     TITLE = {Qualitative aspects of counting real rational curves on real
              {K}3 surfaces},
   JOURNAL = {Geom. Topol.},
  FJOURNAL = {Geometry \& Topology},
    VOLUME = {21},
      YEAR = {2017},
    NUMBER = {1},
     PAGES = {585--601},
      ISSN = {1465-3060,1364-0380},
   MRCLASS = {14J28 (14N35 14P05)},
  MRNUMBER = {3608720},
MRREVIEWER = {Alex\ Degtyarev},
       DOI = {10.2140/gt.2017.21.585},
       URL = {https://doi.org/10.2140/gt.2017.21.585},
}

@article {Welschinger05,
    AUTHOR = {Welschinger, Jean-Yves},
     TITLE = {Invariants of real symplectic 4-manifolds and lower bounds in
              real enumerative geometry},
   JOURNAL = {Invent. Math.},
  FJOURNAL = {Inventiones Mathematicae},
    VOLUME = {162},
      YEAR = {2005},
    NUMBER = {1},
     PAGES = {195--234},
      ISSN = {0020-9910,1432-1297},
   MRCLASS = {53D45 (14N10 14N35 14P25)},
  MRNUMBER = {2198329},
MRREVIEWER = {Eugenii\ Shustin},
       DOI = {10.1007/s00222-005-0445-0},
       URL = {https://doi.org/10.1007/s00222-005-0445-0},
}

@article {Mikhalkin17,
    AUTHOR = {Mikhalkin, Grigory},
     TITLE = {Quantum indices and refined enumeration of real plane curves},
   JOURNAL = {Acta Math.},
  FJOURNAL = {Acta Mathematica},
    VOLUME = {219},
      YEAR = {2017},
    NUMBER = {1},
     PAGES = {135--180},
      ISSN = {0001-5962,1871-2509},
   MRCLASS = {14T05 (14H50)},
  MRNUMBER = {3765660},
MRREVIEWER = {Hannah\ Markwig},
       DOI = {10.4310/ACTA.2017.v219.n1.a5},
       URL = {https://doi.org/10.4310/ACTA.2017.v219.n1.a5},
}

@article {Sottile00,
    AUTHOR = {Sottile, Frank},
     TITLE = {Real rational curves in {G}rassmannians},
   JOURNAL = {J. Amer. Math. Soc.},
  FJOURNAL = {Journal of the American Mathematical Society},
    VOLUME = {13},
      YEAR = {2000},
    NUMBER = {2},
     PAGES = {333--341},
      ISSN = {0894-0347,1088-6834},
   MRCLASS = {14M15 (14N35 14P05)},
  MRNUMBER = {1706484},
MRREVIEWER = {Raquel\ Mallavibarrena},
       DOI = {10.1090/S0894-0347-99-00323-9},
       URL = {https://doi.org/10.1090/S0894-0347-99-00323-9},
}

@article {GKLRSV17,
    AUTHOR = {Gallet, Matteo and Koutschan, Christoph and Li, Zijia and
              Regensburger, Georg and Schicho, Josef and Villamizar, Nelly},
     TITLE = {Planar linkages following a prescribed motion},
   JOURNAL = {Math. Comp.},
  FJOURNAL = {Mathematics of Computation},
    VOLUME = {86},
      YEAR = {2017},
    NUMBER = {303},
     PAGES = {473--506},
      ISSN = {0025-5718,1088-6842},
   MRCLASS = {70B15 (12Y05 14P05 16Z05 68W30 70G55)},
  MRNUMBER = {3557808},
MRREVIEWER = {Peter\ S.\ Donelan},
       DOI = {10.1090/mcom/3120},
       URL = {https://doi.org/10.1090/mcom/3120},
}

@article {LSS19,
    AUTHOR = {Li, Zijia and Schicho, Josef and Schr\"{o}cker, Hans-Peter},
     TITLE = {Factorization of motion polynomials},
   JOURNAL = {J. Symbolic Comput.},
  FJOURNAL = {Journal of Symbolic Computation},
    VOLUME = {92},
      YEAR = {2019},
     PAGES = {190--202},
      ISSN = {0747-7171,1095-855X},
   MRCLASS = {16Z05 (16K99)},
  MRNUMBER = {3907352},
MRREVIEWER = {Michela\ Ceria},
       DOI = {10.1016/j.jsc.2018.02.005},
       URL = {https://doi.org/10.1016/j.jsc.2018.02.005},
}

@article {KM96,
    AUTHOR = {Kapovich, Michael and Millson, John J.},
     TITLE = {The relative deformation theory of representations and flat
              connections and deformations of linkages in constant curvature
              spaces},
   JOURNAL = {Compositio Math.},
  FJOURNAL = {Compositio Mathematica},
    VOLUME = {103},
      YEAR = {1996},
    NUMBER = {3},
     PAGES = {287--317},
      ISSN = {0010-437X,1570-5846},
   MRCLASS = {58H15 (17B55 58D27)},
  MRNUMBER = {1414592},
MRREVIEWER = {Athanase\ Papadopoulos},
       URL = {http://www.numdam.org/item?id=CM_1996__103_3_287_0},
}

@book {Farin97,
    AUTHOR = {Farin, Gerald},
     TITLE = {Curves and surfaces for computer-aided geometric design},
    SERIES = {Computer Science and Scientific Computing},
   EDITION = {Fourth},
      NOTE = {A practical guide,
              Chapter 1 by P. B\'ezier; Chapters 11 and 22 by W. Boehm,
              With 1 IBM-PC floppy disk (3.5 inch; HD)},
 PUBLISHER = {Academic Press, Inc., San Diego, CA},
      YEAR = {1997},
     PAGES = {xviii+429},
      ISBN = {0-12-249054-1},
   MRCLASS = {65D17 (65-01)},
  MRNUMBER = {1412572},
}

@article {MR1575631,
    AUTHOR = {Kempe, A. B.},
     TITLE = {On a {G}eneral {M}ethod of describing {P}lane {C}urves of the
              nth degree by {L}inkwork},
   JOURNAL = {Proc. Lond. Math. Soc.},
  FJOURNAL = {Proceedings of the London Mathematical Society},
    VOLUME = {7},
      YEAR = {1875/76},
     PAGES = {213--216},
      ISSN = {0024-6115},
   MRCLASS = {99-04},
  MRNUMBER = {1575631},
       DOI = {10.1112/plms/s1-7.1.213},
       URL = {https://doi.org/10.1112/plms/s1-7.1.213},
}
\end{document}